\documentclass[12pt]{amsart}
\usepackage{amscd,amsthm,verbatim}
\newcommand{\const}{\operatorname{const.}}

\newcommand{\R}{{\mathbb R}}

\newcommand{\Ric}{\operatorname{Ric}}

\newcommand{\Rm}{\operatorname{Rm}}

\newcommand{\supp}{\operatorname{supp}}

\newcommand{\Z}{{\mathbb Z}}

\numberwithin{equation}{section}
\theoremstyle{plain}

\newtheorem{theorem}[equation]{Theorem}
\newtheorem{proposition}[equation]{Proposition}

\newtheorem{question}[equation]{Question}

\theoremstyle{remark}

\newtheorem{remark}[equation]{Remark}
\newtheorem{example}[equation]{Example}
\usepackage{graphicx}
\input{amssym.def}
\input{amssym}
\input{epsf}
\usepackage{a4wide}

\begin{document}

\title[Backreaction]
      {Backreaction in the future behavior of an expanding
        vacuum spacetime}

\author{John Lott}
\address{Department of Mathematics\\
University of California, Berkeley\\
Berkeley, CA  94720-3840\\
USA} \email{lott@berkeley.edu}

\thanks{Research partially supported by NSF grant
DMS-1510192}
\date{November 29, 2017}

\begin{abstract}
  We perform a rescaling analysis to analyze the future behavior of
  a class of $T^2$-symmetric vacuum spacetimes. We show
  that on the universal cover,
  there is $C^0$-convergence to a spatially homogeneous
  spacetime that does not satisfy the vacuum Einstein equations.
\end{abstract}

\maketitle


\section{Introduction} \label{sec1}

In this paper we study the future behavior of an expanding vacuum
spacetime $(M, g)$ with compact spatial slices.  A basic question is
whether the gravitational dynamics, in the form of the equation
$\Ric(g) = 0$, force the solution to approach a
locally spatially homogeneous spacetime in the future; see
\cite[Part I]{Ringstrom (2013)} for discussion. To make the question precise,
one must say what sort of limit one is considering.

We take the viewpoint that the relevant notion of convergence
is that of a sequence of pointed
vacuum spacetimes. Details are in Section \ref{sec2}
but to give the idea, let $\{p_i\}_{i=1}^\infty$ be a sequence of
points in $M$ going to future infinity. Let $\{c_i\}_{i=1}^\infty$
be a sequence of positive numbers.  Then $\{(M, c_i g, p_i)\}_{i=1}^\infty$ is
a sequence of pointed vacuum spacetimes and we can ask whether
there is a limit $(M_\infty, g_\infty, p_\infty)$ in the pointed
sense.  The latter roughly means that we compare neighborhoods
of $p_i$ of an arbitrary but fixed size, as $i \rightarrow \infty$,
to the corresponding neighborhood of $p_\infty$.
This notion is prevalent in Riemannian geometry and Ricci flow.

One basic issue is that the coordinates used to compute the future asymptotics
of $g$ may not be well adapted to describe the
geometry around $p_i$ for large $i$.  Hence in the definition of
convergence, before taking a limit
one allows $i$-dependent changes of coordinates. One
can think of taking normal coordinates around $p_i$.

There is some freedom in the choice of parameters $\{c_i\}_{i=1}^\infty$,
which determine the scales at which we are making comparisons.
They should have
engineering dimension $time^{-2}$ or $distance^{-2}$, so that
$c_i g_i$ is dimensionless. By a ``type-III rescaling'' we mean
that $c_i$ is constructed using the
proper time of $p_i$ from a fixed hypersurface,
or the Hubble time $t = - \frac{3}{H}$ of $p_i$ with respect to a
constant mean curvature (CMC) spatial foliation with mean curvature
function $H : (T_0, \infty) \rightarrow (H_0, \infty)$, where
$H_0 < 0$. (The negativity of $H$ is the expanding nature of the
spacetime.)

One must also specify the sense in which the metric tensors converge.
In
\cite{Anderson (2001),Lott (2017)}, it was shown that
in the case of a CMC foliation, if the curvature has quadratic decay
in the Hubble time then, after passing to a subsequence,
there is a limit of the metrics in the pointed weak
$W^{2,q}$-topology, for any $q \in [1, \infty)$, with the
  limit being a vacuum spacetime. 
  In such a case, we can also assume that the metrics converge in the pointed
  $C^{1,\alpha}$-topology for any $\alpha \in (0,1)$.

  In this paper we look instead at type-III rescalings of vacuum spacetimes
  that may not satisfy the curvature decay condition. Although a limit
  is no longer guaranteed, we can still ask whether the pointed spacetimes
  have a limit $(M_\infty, g_\infty, p_\infty)$,
  say in the pointed $C^0$-topology. If
  $(M_\infty, g_\infty)$ exists, and is locally spatially homogeneous,
  then it makes sense to say that
  the original $(M,g)$ approaches a locally spatially
  homogeneous spacetime in the
  $C^0$-topology, along the sequence $\{p_i\}$.
  We will actually pass to the universal cover $\widetilde{M}$ and ask
  whether the lifted spacetime $(\widetilde{M}, \widetilde{g})$ approaches
  a spatially homogeneous spacetime.
  It is quite possible
  that $(\widetilde{M}, \widetilde{g})$
  approaches a spatially homogeneous spacetime
  in the $C^0$-topology, but
  not in some stronger topology.

  We perform this rescaling analysis for a class of vacuum spacetimes
  with compact spatial slices diffeomorphic to $T^3$, and invariance
  under the action of the group $T^2$. Such a spacetime
  is {\em polarized} if the Killing fields can be taken to be
  orthogonal, i.e. if the $T^2$-invariance can be promoted to an
  $(O(2) \times O(2))$-invariance. A $T^2$-symmetric spacetime
  has a {\em twist constant} $K$; if $K=0$ then the spacetime is
  {\em Gowdy}. Future asymptotics of $T^2$-invariant spacetimes were
  considered
  by Ringstr\"om in the Gowdy case \cite{Ringstrom (2004),Ringstrom (2006)} and
  the nonGowdy case \cite{Ringstrom (2015)}. More precise
  asymptotics were obtained by LeFloch and Smulevici in the
  polarized nonGowdy case, for initial data that is sufficiently close
  to the asymptotic regime \cite{LeFloch-Smulevici (2016)}.

  \begin{proposition} \label{prop1.1}
    After passing to the universal cover,
    any vacuum spacetime of the type considered in
    \cite{LeFloch-Smulevici (2016)} has
    a smooth type-III rescaling limit $g_\infty$, in the pointed
    $C^0$-topology.  The Lorentzian metric $g_\infty$
    is spatially homogeneous but does not satisfy
    the vacuum Einstein equations.
    \end{proposition}
  
  The effective stress-energy tensor $T = \Ric_{g_\infty}
  - \frac12 R_{g_\infty} g_\infty$ of $g_\infty$ vanishes except
  for the $T_{00}$ component, which is positive.
  (We do not claim that $T$ has any physical meaning.)
  The fact that $g_\infty$ does not
  satisfy the vacuum Einstein equations implies that the convergence
  cannot be in the pointed $(C^0 \cap H^1)$-topology, as otherwise
  the vacuum Einstein equations would make sense weakly and pass
  to the limit; c.f. \cite{Lott (2016)}.

  That a limit of vacuum spacetimes can
  have a nonzero stress-energy tensor is called
  backreaction
  \cite{Buchert-Rasanen (2012),Green-Wald (2011),Huneau-Luk (2017)}.
  In effect, fluctations of
  the geometry, with increasing frequency, can average out to zero in
  some parts of the Einstein equations, but give a nonzero
  contribution through nonlinearities to other parts.
  This phenomenon of increasing
  fluctuations also arose in the analysis of expanding spacetimes
  \cite{Lott (2017),Ringstrom (2004),Ringstrom (2006)}.
  In \cite{Green-Wald (2011)}, a framework was developed to analyze
  backreaction, with one of the main conclusions being that 
  the effective
  stress-energy tensor is trace-free.  We see that the framework
  of \cite{Green-Wald (2011)} does not apply to our rescaling examples.
  
  The fact that $T_{00}$ is positive in Proposition \ref{prop1.1}
  makes one wonder how generally a limiting stress-energy tensor
  satisfies some positive energy condition.
    Motivated by the results of \cite{Bamler (2016), Gromov (2014)},
    in Section \ref{sec4} we raise a purely Riemannian
    question about the behavior of scalar curvature
  when taking a $C^0$-limit of Riemannian metrics.
  In Proposition \ref{prop4.3}
  we show that a positive answer to this question implies that
  if a sequence of CMC vacuum spacetimes converges in the
  pointed weak $H^1$-topology and the pointed $C^0$-topology, then
  the limiting spacetime has a nonnegative energy density.

  The structure of the paper is the following.  In Section \ref{sec2}
  we discuss rescaling limits of expanding vacuum spacetimes.
  In Section \ref{sec3} we analyze the polarized $T^2$-symmetric
  spacetimes of \cite{LeFloch-Smulevici (2016)}. Section \ref{sec4}
  has the link to questions of scalar curvature in Riemannian
  geometry. Section \ref{sec5} has a short discussion of the results
  of the paper.

  I thank Hans Ringstr\"om for references to the literature, and
  the referees for useful comments.
  
  \section{Rescaling limits} \label{sec2}

  In this section we discuss notions of pointed convergence
  and rescaling for
  spacetimes.  These notions are not new, at least
  in spirit; c.f. \cite{Anderson (2001),Anderson (2003)}.
  
\subsection{Rescaling limits of spacetimes} \label{subsec2.1}

Let $\{(M_i,g_i)\}_{i=1}^\infty$ be a sequence of $(n+1)$-dimensional
Lorentzian manifolds. For the moment, we do not specify the regularity of
the metrics.
Let $p_i \in M_i$ be a basepoint.
Let $(M_\infty, g_\infty)$ be another such Lorentzian
manifold, with basepoint $p_\infty \in M_\infty$. We say that
$\lim_{i \rightarrow \infty} (M_i, g_i, p_i) = (M_\infty, g_\infty, p_\infty)$
if there is
\begin{itemize}
\item An exhaustion of $M_\infty$ by compact codimension-zero
submanifolds-with-boundary $p_\infty \in K_1 \subset K_2 \subset \ldots$,
and \\
\item Maps $\phi_{i,j} : K_j \rightarrow M_i$,
with $\phi_{i,j}(p_\infty) = p_i$, that are diffeomorphisms onto their
images for large $i$, such that \\
\item For all $j$, we have
$\lim_{i \rightarrow \infty} \phi_{i,j}^* g_i = g_\infty$ on $K_j$.
\end{itemize}
Here the notion of convergence of metrics depends on the topology
that we want to consider, e.g. $C^0$, $C^k$, $C^\infty$, $W^{q,k}$,
etc.
If each $(M_i, g_i)$ is a $C^2$-regular vacuum spacetime,
i.e. $\Ric(g_i) = 0$, and $g_\infty$ is $C^2$-regular,
then $(M_\infty, g_\infty)$ is a vacuum
spacetime provided that the metric convergence is $C^0 \cap H^1$,
since the Ricci-flat condition then makes sense weakly.

If we start with a single smooth Lorentzian manifold $(M,g)$, and a
sequence $\{p_i\}_{i=1}^\infty$ in $M$, then we may want to take
$M_i = M$ and $g_i = c_i g$ for some constants $c_i > 0$. The goal is to find
constants $c_i$ and maps $\phi_{i,j} : K_j \rightarrow M$ so that
there is a limit $(M_\infty, g_\infty)$. We want $c_i$ to have
engineering dimension $time^{-2}$ or $distance^{-2}$, so that
$g_i$ is dimensionless.
If one takes
$\{c_i\}_{i=1}^\infty$ increasing
sufficiently quickly then one can always get a flat limit in the smooth
topology, but this would
be considered uninteresting.

\subsection{Rescaling limits of CMC spacetimes} \label{subsec2.2}

Going back to the sequence $\{(M_i, g_i)\}_{i=1}^\infty$,
suppose that each $(M_i, g_i)$ is smooth and
that there is a globally hyperbolic
foliation $M_i = (T_i, \infty) \times X_i$ by
spatial hypersurfaces.
We can generally perform spatial diffeomorphisms to write
  \begin{equation} \label{2.1}
g_i = - L_i^2(t) dt^2 + h_i(t),
  \end{equation}
  where $L_i(t) \in C^\infty(X_i)$ is the lapse function
  and $h_i(t)$ is a Riemannian metric on
  $X_i$; for example, we can always put $g_i$ in this form if $X_i$ is compact.
 Write
 $p_i = (t_i, x_i)$.
 
 Let
 $M_\infty = I_\infty \times Y$ be a putative limit space, for some
 open interval $I_\infty \subset \R$,
 with basepoint $(u_\infty, y_\infty)$.
  One may
  wish to consider comparison maps that preserve the structure (\ref{2.1}).
  To do so, let $u_\infty \in C_1 \subset C_2 \ldots$ be an exhaustion
  of $I_\infty$ by compact intervals.  Let $\sigma_{i,j} : C_j \rightarrow
 (T_i, \infty)$ be a map with $\sigma_{i,j}(u_\infty) = t_i$
 that is a diffeomorphism to its image for large $i$.
 Let $y_\infty \in K^\prime_1 \subset
  K^\prime_2 \subset \ldots$ be an exhaustion of $Y$ by
  compact codimension-zero submanifolds-with-boundary.
  Put $K_j = C_j \times K^\prime_j$.
  Given maps $\eta_{i,j} : K^\prime_j \rightarrow X_i$
  with $\eta_{i,j}(y_\infty) = x_i$ that
  are diffeomorphisms to their images for large $i$, we define
  $\phi_{i,j} : K_j \rightarrow M_i$ by
  $\phi_{i,j}(u,y) = (\sigma_{i,j}(u), \eta_{i,j}(y))$.
Then $\phi_{i,j}^* g_i$ has the form (\ref{2.1}).
If the convergence $\lim_{i \rightarrow \infty} \phi_{i,j}^* g_i =
g_\infty \big|_{K_j}$ is $C^0$ for each $j$, then $g_\infty$ also has the form
(\ref{2.1}).

A special case is when for each $t \in (T_i, \infty)$, the hypersurface
$\{t\} \times X_i$ has constant mean curvature $H_i(t)$.
Then $\phi_{i,j}^* g_i$ also has a constant mean curvature (CMC)
foliation.
If the convergence $\lim_{i \rightarrow \infty} \phi_{i,j}^* g_i =
g_\infty \big|_{K_j}$ is $C^1$ for each $j$
then $g_\infty$ acquires a CMC foliation.

We say that $(M, g)$ is an expanding CMC spacetime if
$H : (T_0, \infty) \rightarrow (H_0, 0)$ is bijective
and increasing, where $H_0 < \infty$.
Then we can assume that the time parameter $t$ of $M$ satisfies
  $t = - \frac{n}{H}$, i.e. $t$ is the Hubble time of the
  CMC foliation.
  
\subsection{Type-III rescalings} \label{subsec2.3}

Continuing with an expanding CMC spacetime $M = (T_0, \infty) \times X$,
parametrized by Hubble time $t$,
and a sequence $p_i = (t_i, x_i)$ with $\lim_{i \rightarrow \infty} t_i =
\infty$,
we say that a type-III rescaling is when $M_i = M$, $X_i = X$,
  $I_\infty = (0, \infty)$, $u_\infty = 1$,
  $C_j = \left[ \frac{1}{j}, j \right]$,
  $\sigma_{i,j}(u) = t_i u$ (for all $i$ sufficiently large that
  $t_i > T_0 j$) and $c_i = t_i^{-2}$. Then
 $u$ is the Hubble time for the CMC foliation of $\phi_{i,j}^* g$.
If the convergence $\lim_{i \rightarrow \infty} \phi_{i,j}^* g =
g_\infty \big|_{K_j}$ is $C^1$ then $u$ is also the Hubble time
for the CMC foliation of $g_\infty$.

Define a curvature norm for $(M, g)$, at a point $m \in M$,
as follows
\cite[(0.7)]{Anderson (2001)}.  Let $\{e_i\}_{i=0}^n$
be an orthonormal basis for $T_mM$ with $e_0 =
(- g(\partial_t, \partial_t))^{- \: \frac12} \partial_t$. Put
\begin{equation} \label{2.2}
|\Rm|(m) = \sqrt{\sum_{\alpha, \beta, \gamma, \delta = 0}^n
  R(e_\alpha, e_\beta, e_\gamma, e_\delta)^2}.
\end{equation}
Assuming that each $(X,h(t))$ is complete,
and $|\Rm| = O \left(t^{-2} \right)$, we can use the
type-III scaling to extract a
subsequential limit $(M_\infty, g_\infty, p_\infty)$
\cite[Corollary 3.4]{Lott (2017)}. The limit is in the
pointed weak $W^{2,q}$-topology for all $q \in [1, \infty)$,
but  $Y$ may be an \'etale groupoid rather than a manifold,
  if there is ``collapsing''.
  In some cases, such as if $X$ is compact and aspherical,
  one can stay in the world of manifolds by lifting $g$ to
  the universal cover $(T_0, \infty)
  \times \widetilde{X}$ and taking a pointed limit there.

  \begin{example} \label{ex2.3}
    Consider a Kasner spacetime $(0, \infty) \times T^n$ with
    metric
\begin{equation} \label{2.4}
g \: = \: - \: \frac{1}{n^2} dt^2 +
\sum_{k=1}^n t^{2p_k} (dx^k)^2,
\end{equation}
where
    $\sum_{k=1}^n p_k = \sum_{k=1}^n p_k^2 = 1$.
    Then $t$ is the Hubble time. Passing to the universal cover,
    we take $x^k \in \R$.
    Put $M_\infty = (0, \infty) \times \R^n$, with
    $p_\infty = (1,0)$. Writing a point $x \in \R^n$ as
    $x = (x^k)$, define $\eta_{i,j}(y) \in \R^n$ to be the
    point whose $k^{th}$-coordinate is
    $t_i^{1-p_k} y^k + x_i^k$. Then $\phi_{i,j}^* g_i$ is the
    Kasner metric on $M_\infty$, now with time parameter $u$.
    Hence in this case the rescaling limit exists on the
    universal cover.
    \end{example}
  
  If the foliation $M = (T_0, \infty) \times X$ may not be
  a CMC foliation, an
  alternative type-III rescaling uses the proper time from
  a fixed hypersurface
\cite{Anderson (2001)}.
Fixing $T_1 \in (T_0, \infty)$, define
  $\tau : [T_1, \infty) \rightarrow \R$ by saying that
    $\tau(t)$ is the maximal length of 
  causal curves from the time-$T_1$ hypersurface to the
  time-$t$ hypersurface. We reparametrize $M$ by $\tau$.
    Putting $\tau_i = \tau(t_i)$, we can rescale using
    $\sigma_{i,j}(u) = \tau_i u$ and $c_i = \tau_i^{-2}$.
  
\subsection{Type-II rescalings} \label{subsec2.4}

Again in the case of an expanding CMC spacetime, parametrized by
Hubble time $t$,
if $|\Rm|$ is not $O \left(t^{-2} \right)$ then it makes sense to do a
  type-II rescaling.
Put $I_\infty = \R$, $u_\infty = 0$ and
$C_j = \left[ -j, j \right]$.
Given a sequence $p_i = (t_i, x_i)$ with $\lim_{i \rightarrow \infty}
t_i = \infty$, put
  $\sigma_{i,j}(u) = |\Rm(p_i)|^{- \: \frac12} u + t_i$
  (for all $i$ sufficiently large that
  $t_i - |\Rm(p_i)|^{- \: \frac12} j > T_0 $) and
$c_i = |\Rm(p_i)|$.
By construction, the curvature tensor of
$g_i$ at $p_i$ has norm one.
With the right choice of $\{p_i\}_{i=1}^\infty$,
there is a subsequential limit
$g_\infty$ in the pointed weak
$W^{2,q}$-topology for any $q \in [1, \infty)$
  \cite[Proposition 2.51]{Lott (2017)}.
    Two caveats must be made. First, in the collapsing case,
    $Y$ may be an \'etale groupoid rather than a manifold.
    Second, the limiting lapse function $L_\infty$ may vanish.
    If this happens then $g_\infty$ is a static solution of the
    constraint equations.

    If the second fundamental form $K$ of
    $g$ satisfies an inequality $|K|^2 \le \const H^2$ then
    $L_\infty > 0$ \cite[Proposition 4.1]{Lott (2017)}.
    If in addition $n=3$ then $g_\infty$ turns out to be a flat static
    spacetime
    \cite[Corollary 4.6]{Lott (2017)}. The interpretation is that
    there are increasing fluctuations of the curvature tensor, at
    least in neighborhoods of the points $p_i$, that average out
    the normalized curvature to become zero in the weak limit.

    \section{Polarized $T^2$-symmetric nonGowdy spacetimes} \label{sec3}
We now take $n=3$ and
    $X = T^3$, with linear coordinates $(\theta, x, y) \in
    (\R/2\pi\Z)^3$. 
    We assume
    that there is an $(O(2) \times O(2))$-symmetry, acting on the
    $(x,y)$-factor.
    Take the time parameter $R$ so that the area of the
    $T^2$-orbit is $R^2$. As in
    \cite[(2-2)]{LeFloch-Smulevici (2016)},
    the metric can be written
    \begin{equation} \label{3.1}
      g = e^{2(\eta - U)} (- \: dR^2 + a^{-2} d\theta^2) +
        e^{2U} (dx + G d\theta)^2 +
        e^{-2U} R^2 (dy + H d\theta)^2,
    \end{equation}
    where $\eta, U, a, G$ and $H$ are functions of $R$ and $\theta$.
    Let $K$ be the twist constant.  
    We assume that $K \neq 0$.
Let $\langle \eta \rangle (R)$ denote the average value of
$\eta(R, \theta)$ with respect to $\theta \in \R/2 \pi \Z$, and
similarly for $\langle U \rangle (R)$.
    From
    \cite[Theorem 7.1 and (2-8)]{LeFloch-Smulevici (2016)},
    if the initial data are close enough to the asymptotic regime then
    the leading asymptotics
    of the metric parameters are 
    \begin{align} \label{3.2}
       |K^2  e^{2\eta} - R^2| = & \: O \left( R^{\frac74} \right), \\
|\eta - \langle \eta \rangle| = & O \left( R^{- \: \frac12} \right), \notag \\
       |a^{-1} - \frac{2}{\sqrt{5}} C_\infty^{\frac12} R^{\frac12}
       {\mathcal L}(\theta)| = &
       O \left( R^{-1} \right), \notag \\
       |U - C_U|  = & \: O \left( R^{- \: \frac12} \right),  \notag \\
       |G - {\mathcal G}(\theta)| = & \: 0, \notag \\
      |H - \frac{4}{K \sqrt{5}} C_\infty^{\frac12} R^{\frac12}
      {\mathcal L}(\theta)| = & \: O \left( R^{\frac14} \right),  \notag 
      \end{align}
where $C_U, C_\infty$ are constants with $C_\infty > 0$, and 
${\mathcal G}(\theta), {\mathcal L}(\theta)$ are functions with
${\mathcal L}(\theta) > 0$.

The constant-$R$ slices are generally not CMC.
The maximal
length of causal curves
between a constant $R_0$-slice and a constant $R$-slice is
asymptotic to
\begin{equation} \label{3.3}
  \int_{R_0}^R e^{\langle \eta \rangle(r) - \langle U \rangle(r)} \: dr \sim
  \const \int_{R_0}^R r \: dr \sim \const R^2;
\end{equation}
c.f. \cite[Pf. of Proposition 3]{Ringstrom (2006)}.
Changing variable from $R$ to $t = R^2$, the asymptotic behavior of $g$ is
\begin{align} \label{3.4}
  g \sim & 
  - \: \frac14 K^{-2} e^{- 2C_U} dt^2 + \frac45 
  K^{- \: 2} e^{-\: 2 C_U}
  C_\infty {\mathcal L}^2 t^{\frac32} d\theta^2 +
  e^{2C_U} (dx + {\mathcal G} d\theta)^2 + \\
  & e^{- \: 2 C_U} t
  \left( dy +
  \frac{4}{K \sqrt{5}} C_\infty^{\frac12}
  {\mathcal L} t^{\frac14} d\theta \right)^2. \notag
  \end{align}
That is, when restricted to a future time interval
$t \in [c, \infty)$, the two sides of (\ref{3.4}) are
$(1 + o(c))$-biLipschitz.

  Let $p_i = (t_i, x_i)$ be a sequence with $\lim_{i \rightarrow \infty}
  t_i = \infty$.  The choice of points $x_i \in T^3$ will be irrelevant.
  Putting $t = t_i u$ gives
\begin{align} \label{3.5}
  t_i^{-2} g \sim & 
  - \: \frac14 K^{-2} e^{- 2C_U} du^2 + \frac45 
  K^{- \: 2} e^{-\: 2 C_U}
  C_\infty {\mathcal L}^2 t_i^{- \: \frac12} u^{\frac32} d\theta^2 +
  e^{2C_U} \left( t_i^{-1} dx + t_i^{-1} {\mathcal G} d\theta \right)^2 + \\
  & e^{- \: 2 C_U} u
  \left(  t_i^{- \: \frac12} dy +
  \frac{4}{K \sqrt{5}} C_\infty^{\frac12}
  {\mathcal L} t_i^{- \: \frac14} u^{\frac14} d\theta \right)^2. \notag
  \end{align}
Passing to the universal cover, we take
$(\theta, x, y) \in \R^3$. For simplicity, we just take the spatial
basepoint to be $0 \in \R^3$.
We define $\widehat{\theta}$, $\widehat{x}$ and $\widehat{y}$ by
$d\widehat{\theta} = t_i^{- \: \frac14} {\mathcal L} d\theta$,
$\widehat{x} = t_i^{-1} x$ and $\widehat{y} = t_i^{- \: \frac12} y$. Then
\begin{align} \label{3.6}
  t_i^{-2} g \sim & 
  - \: \frac14 K^{-2} e^{- 2C_U} du^2 + \frac45 
  K^{- \: 2} e^{-\: 2 C_U}
  C_\infty u^{\frac32} d\widehat{\theta}^2 +
  e^{2C_U} (d\widehat{x} + {\mathcal G} {\mathcal L}^{-1}
  t_i^{- \frac34} d\widehat{\theta})^2 + \\
  & e^{- \: 2 C_U} u
  \left(  d\widehat{y} +
  \frac{4}{K \sqrt{5}} C_\infty^{\frac12}
  u^{\frac14} d\widehat{\theta} \right)^2. \notag
  \end{align}
Since $\lim_{i \rightarrow \infty} t_i = \infty$, there is a
limit $\lim_{i \rightarrow \infty} t_i^{-2} g_i = g_\infty$
in the pointed $C^0$-topology:
\begin{align} \label{3.7}
  g_\infty = & 
  - \: \frac14 K^{-2} e^{- 2C_U} du^2 + \frac45 
  K^{- \: 2} e^{-\: 2 C_U}
  C_\infty u^{\frac32} d\widehat{\theta}^2 +
  e^{2C_U} d\widehat{x}^2 + \\
  & e^{- \: 2 C_U} u
  \left(  d\widehat{y} +
  \frac{4}{K \sqrt{5}} C_\infty^{\frac12}
  u^{\frac14} d\widehat{\theta} \right)^2. \notag
\end{align}
We see that $g_\infty$ has a spatial $\R^3$-symmetry.
Redefining 
$u = \widehat{R}^2$
gives
\begin{equation} \label{3.8}
  g_\infty = e^{2(\widehat{\eta} - \widehat{U})} (- \: d\widehat{R}^2 +
  \widehat{a}^{-2}
  d\widehat{\theta}^2) +
        e^{2\widehat{U}} (d\widehat{x} + \widehat{G} d\widehat{\theta})^2 +
        e^{-2\widehat{U}}
        \widehat{R}^2 (d\widehat{y} + \widehat{H} d\widehat{\theta})^2,
\end{equation}
where
\begin{align} \label{3.9}
       e^{2\widehat{\eta}} = & \: K^{-2} \widehat{R}^2, \\
       \widehat{a}^{-1} = &  \frac{2}{\sqrt{5}}
       C_\infty^{\frac12} \widehat{R}^{\frac12},
       \notag \\
       \widehat{U} = & C_U,  \notag \\
       \widehat{G} = & \: 0, \notag \\
       \widehat{H} = & \frac{4}{K \sqrt{5}} C_\infty^{\frac12}
       \widehat{R}^{\frac12}.  \notag 
\end{align}
To check whether $g_\infty$ satisfies the vacuum Einstein equations,
we can plug (\ref{3.9}) into \cite[(2-3)-(2-8)]{LeFloch-Smulevici (2016)}.
One finds that these equations are satisfied except for the constraint
equation
\cite[(2-6)]{LeFloch-Smulevici (2016)}, which instead becomes
\begin{equation} \label{3.10}
  \widehat{\eta}_{\widehat{R}} + \frac{K^2}{4 \widehat{R}^3}
  e^{2 \widehat{\eta}} -
  \widehat{a} \widehat{R}
  \left( \widehat{a}^{-1} \widehat{U}_{\widehat{R}}^2 +
  \widehat{a} \widehat{U}_{\widehat{\theta}}^2 \right) =
  \frac{5}{4\widehat{R}}.
  \end{equation}
The left-hand side of (\ref{3.10}) is proportionate to
$(\Ric_\infty - \frac12 R_\infty g_\infty)\left( \partial_{\widehat{R}},
\partial_{\widehat{R}} \right)$.
We conclude that $g_\infty$ satisfies the Einstein equations,
except for the nonvanishing of
$(\Ric_\infty - \frac12 R_\infty g_\infty)\left( \partial_{\widehat{R}},
\partial_{\widehat{R}} \right)$.

\begin{remark} \label{rmk3.11}
  One can also consider rescaling limits of Gowdy spacetimes, i.e.
  $T^2$-symmetric spacetimes with spatial slices
  diffeomorphic to $T^3$, and vanishing twist constant.
  For such spacetimes, the curvature decays like the inverse
  square of the proper time function, as measured from a fixed
  hypersurface \cite[Theorem 2]{Ringstrom (2006)}. Hence we
  expect that they have type-III rescaling limits that are
  vacuum spacetimes.  If
  the metric is independent of the parameter $\theta$ of
  $S^1 = T^3/T^2$ then before rescaling,
  the solution on the universal cover is a
  spatially homogeneous Kasner spacetime.
  In this case the rescaling limit exists and is also a Kasner spacetime;
  see Example \ref{ex2.3}.
  If the metric is not $\theta$-independent then
  asymptotics were given in \cite{Ringstrom (2004),Ringstrom (2006)}.
  Some rough calculations indicate that
  the rescaled metrics should approach a flat metric (in the
  {\em weak} $W^{2,q}$-topology).  However, it does not seem to be possible
  to prove this rigorously from the known asymptotics.

  Another interesting vacuum spacetime is the Bianchi VIII solution. The
  so-called non-NUT type does not have curvature that decays like
  the inverse square of the Hubble time (or the proper time)
  \cite[Theorem 3]{Ringstrom (2006)}.
    Based on some calculations in terms of
  the coordinates from \cite[Section 4.3]{Glickenstein (2008)} or
  \cite[Section 3.3.5]{Lott (2007)},
  there does not appear to
  be a type-III rescaling limit in the $C^0$-topology.
  \end{remark}

\section{Nonnegativity of induced energy density} \label{sec4}

We recall that in Subsection \ref{subsec2.3}, there was a subsequential
rescaling limit of a
CMC vacuum spacetime with quadratic curvature decay,
that exists in the pointed weak $W^{2,q}$-topology
for all $q \in [1, \infty)$. This is related to the fact from Riemannian
  geometry that a sequence of complete pointed Riemannian manifolds,
  with uniformly bounded curvature, has a subsequential limit in the
  pointed weak $W^{2,q}$-topology.

  In Riemannian geometry, if one weakens the curvature assumptions to
  a uniform lower bound on the Ricci curvature, and a uniform lower
  bound on the injectivity radius, then there is 
a subsequential limit in the pointed weak
$W^{1,q}$-topology for all $q \in [1, \infty)$, and hence also a
  subsequential limit in the pointed $C^{\alpha}$-topology for all
    $\alpha \in (0,1)$ \cite{Anderson-Cheeger (1992)}.
    Motivated by this, we consider a sequence of 
    CMC vacuum spacetimes $(T_i, \infty) \times X_i$, as in Subsection \ref{subsec2.2},
    so that for each $u \in I_\infty$,
\begin{itemize}
\item The pullback metrics $h_i(u)$
  and the pullback lapse functions
  $L_i(u)$ converge in the pointed weak $H^1$-topology
    and
    the pointed $C^0$-topology, and
    \item The pullback
      second fundamental forms
      $K_i(u)$ converge in the pointed weak
      $L^2$-topology.
\end{itemize}

    In Riemannian geometry, there is a general principle that curvature can
    only go up when taking limits.  In the case of scalar curvature, a
    precise statement along these lines is the following result.

\begin{theorem} \label{thm4.1} \cite{Bamler (2016),Gromov (2014)}
    Let $Y$ be a smooth manifold. Given $\kappa \in C(Y)$, let
    $\{{\frak g}_i\}_{i=1}^\infty$ be a sequence of $C^2$-regular
    Riemannian metrics on $Y$ with scalar curvature function
    bounded below by $\kappa$. If $\{{\frak g}_i\}_{i=1}^\infty$
    converges on compact subsets in the $C^0$-topology
    to a $C^2$-regular Riemannian metric ${\frak g}_\infty$,
       then ${\frak g}_\infty$ has scalar curvature function
    bounded below by $\kappa$.
    \end{theorem}

    \begin{question} \label{ques4.2}
      Let $Y$ be a smooth manifold. Let
    $\{{\frak g}_i\}_{i=1}^\infty$ be a sequence of $C^2$-regular
    Riemannian metrics on $Y$ that $C^0$-converges on compact subsets
    to a $C^2$-regular Riemannian metric ${\frak g}_\infty$ on $Y$.
    Suppose that on any compact subset, the scalar curvatures of
    $\{{\frak g}_i\}_{i=1}^\infty$ are uniformly bounded below.
    Is it true that for every nonnegative compactly supported smooth
    density $\omega$ on $Y$, the scalar
    curvature $R_{{\frak g}_\infty}$ satisfies
    $\int_Y R_{{\frak g}_\infty} \omega \ge
    \liminf_{i \rightarrow \infty}
    \int_X R_{{\frak g}_i} \omega$?
      \end{question}

    A positive answer to Question \ref{ques4.2} clearly implies Theorem
    \ref{thm4.1}.

    \begin{proposition} \label{prop4.3}
      Consider a sequence of expanding 
    CMC vacuum spacetimes $(T_i, \infty) \times X_i$, as in Subsection \ref{subsec2.2},
    that by assumption converges in the sense of the bulletpoints above
to a
    CMC spacetime $I_\infty \times Y$,
    equipped with a $C^2$-regular metric $g_\infty$
    that is also parametrized by the Hubble time
  Then a positive answer to Question \ref{ques4.2}
  implies that $\left( \Ric_{g_\infty} -
  \frac12 R_{g_\infty} g_\infty \right)(\partial_u, \partial_u) \ge 0$. 
      \end{proposition}
    \begin{proof}
On a given $u$-slice of the limit space, 
the Gauss equation gives
      \begin{equation} \label{4.4}
        ( \Ric_{g_\infty} -
        \frac12 R_{g_\infty} g_\infty )(\partial_u, \partial_u) =
\frac12 \left( R_{h_\infty} - |K_\infty|^2 + H_\infty^2 \right).
      \end{equation}
      By assumption, $H_\infty = - \frac{n}{u}$. Choose a
      nonnegative compactly supported smooth density $\omega$
      on $Y$.
      Then
      \begin{equation} \label{4.5}
      \int_Y \left( (\Ric_{g_\infty} -
      \frac12 R_{g_\infty} g_\infty)(\partial_u, \partial_u) \right)
      \omega =
      \frac12
      \int_Y \left( R_{h_\infty} - |K_\infty|^2 + \frac{n^2}{u^2}  \right)
      \omega.
      \end{equation}
      For large $i$, we can pullback
      $h_i$ and $K_i$ to $\supp(\omega)$, so we assume that
      everything lives on $Y$.
      The constraint equations give
      \begin{equation} \label{4.6}
        R_{h_i} - |K_i|^2 + \frac{n^2}{u^2} =
        R_{h_i} - |K_i|^2 + H_i^2 = 0.
        \end{equation}
In particular, $R_{h_i}$ is bounded below in terms of $u$.
      A positive answer to Question \ref{ques4.2} implies that
      \begin{equation} \label{4.7}
        \int_Y   R_{h_\infty} 
        \omega \ge \liminf_{i \rightarrow \infty}
        \int_Y   R_{h_i} \omega =
        \liminf_{i \rightarrow \infty}
        \int_Y   \left( |K_i|^2 - \frac{n^2}{u^2} \right) \omega.
      \end{equation}
      Then
      \begin{align} \label{4.8}
& \int_Y \left(  (\Ric_{g_\infty} -
      \frac12 R_{g_\infty} g_\infty)(\partial_u, \partial_u) \right)
      \omega \ge \\
      & \frac12 \liminf_{i \rightarrow \infty} \int_Y
      \left( |K_i|^2 - |K_\infty|^2 \right) \omega = \notag \\
      & \frac12 \liminf_{i \rightarrow \infty} \int_Y
      \left( |K_i - K_\infty|^2 + 2
      \langle K_i - K_\infty, K_{\infty} \rangle \right)
      \omega. \notag
             \end{align}
    As $\lim_{i \rightarrow \infty} K_i = K_\infty$ in the weak $L^2$-topology
    on $\supp(\omega)$, it follows that
    \begin{equation} \label{4.9}
      \lim_{i \rightarrow \infty}
      \int_Y \langle K_i - K_\infty, K_{\infty} \rangle \omega = 0.
      \end{equation}
    From (\ref{4.8}), we obtain that
    \begin{equation} \label{4.10}
      \int_Y \left(  (\Ric_{g_\infty} -
      \frac12 R_{g_\infty} g_\infty)(\partial_u, \partial_u) \right)
      \omega \ge 0
    \end{equation}
    for every nonnegative compactly supported smooth density $\omega$
    on $Y$. This implies that
    $(\Ric_{g_\infty} -
      \frac12 R_{g_\infty} g_\infty)(\partial_u, \partial_u) \ge 0$.
    \end{proof}

    \section{Discussion} \label{sec5}

    In this paper we described a notion of rescaling limits for
    Lorentzian spacetimes. For a class of $T^2$-symmetric
    vacuum spacetimes,
    we showed that on the universal cover, there is a rescaling limit
    in the pointed $C^0$-topology that is smooth and spatially
    homogeneous, but does not satisfy the vacuum Einstein equations.

    The paper \cite{Green-Wald (2011)} showed that under certain
    assumptions, a weak limit of a $1$-parameter family of vacuum
    spacetimes has an effective stress-energy tensor that is
    traceless.  The assumptions are about the asymptotics
    of the metric tensors as the parameter $\lambda$ goes to zero;
    we refer to \cite{Green-Wald (2011)} for the details.
    In our examples, the effective stress-energy
    tensor is not traceless. Hence the assumed asymptotics of
    \cite{Green-Wald (2011)} do not hold.  One could try to
    perform a more detailed analysis.

    More generally, one could look at rescaling limits of other
    solutions of the Einstein equations.  In this paper, we
    focused on the future behavior of expanding solutions. One
    could also consider rescaling limits as a singularity develops
    or, similarly, as one goes backward in time toward an initial
    singularity. There is some relation here to the paper
    \cite{Andersson-vanElst-Lim-Uggla (2005)}.


\begin{thebibliography}{10}

    \bibitem{Anderson (2001)}
      M. Anderson, ``On long-time evolution in general relativity and geometrization of 3-manifolds'',
      Comm. Math. Phys. 222, p. 533-567 (2001)

    \bibitem{Anderson (2003)} M. Anderson, ``Regularity for Lorentz metrics under curvature bounds'',
      J. Math. Phys 44, p. 2994-3012 (2003) 
      
    \bibitem{Anderson-Cheeger (1992)} M. Anderson and J. Cheeger,
      ``$C^\alpha$-compactness for manifolds with Ricci curvature and
      injectivity radius bounded below'', J. Diff. Geom. 35, p. 265-281 (1992)

    \bibitem{Andersson-vanElst-Lim-Uggla (2005)}
      L. Andersson, H. van Elst, W. Lim and C. Uggla,
      ``Asymptotic silence of generic cosmological singularities'',
      Phys. Rev. Lett. 94, 051101 (2005)
      
    \bibitem{Bamler (2016)} R. Bamler, ``A Ricci flow proof of a result by Gromov on lower bounds for scalar curvature'', Math. Res. Lett. 23, p. 325-337
      (2016)

    \bibitem{Buchert-Rasanen (2012)}
      T. Buchert and S. R\"as\"anen, ``Backreaction in late time cosmology'',
      Annual Review of Nucl. and Part. Science 62, p. 57-79 (2012)

    \bibitem{Glickenstein (2008)} D. Glickenstein,
      ``Riemannian groupoids and solitons for three-dimensional homogeneous Ricci and cross curvature flows'',
Int. Math. Res. Not. IMRN 2008, no. 12, Article ID rnn034, 49 pages (2008)

\bibitem{Green-Wald (2011)} S. Green and R. Wald,
``A new framework for analyzing the effects of small scale inhomogeneities in cosmology'',
  Phys. Rev. D83, 084020 (2011)

    \bibitem{Gromov (2014)}
      M. Gromov, ``Dirac and Plateau billiards in domains with corners'',
      Cent. Eur. J. Math. 12, p. 1109-1156 (2014)

    \bibitem{Huneau-Luk (2017)} C. Huneau and J. Luk,
      ``High-frequency backreaction for the Einstein equations under
      polarized $U(1)$ symmetry'',
      preprint, https://arxiv.org/abs/1706.09501 (2017)
      
    \bibitem{LeFloch-Smulevici (2016)} P. LeFloch and J. Smulevici,
      ``Future asymptotics and geodesic completeness of polarized
      $T^2$-symmetric spacetimes'', Analysis \& PDE 9, p. 363-395 (2016) 

    \bibitem{Lott (2007)} J. Lott,
      ``On the long-time behavior of type-III Ricci flow solutions'',
      Math. Annalen 339, p. 627-666  (2007)

    \bibitem{Lott (2016)} J. Lott,
      ``Ricci measure for some singular Riemannian metrics'',
      Math. Annalen 365, p. 449-471 (2016) 
      
\bibitem{Lott (2017)} J. Lott, ``Collapsing in the Einstein flow'', preprint,
  https://arxiv.org/abs/1701.05150 (2017)

\bibitem{Ringstrom (2004)} H. Ringstr\"om,
  ``On a wave map equation arising in general relativity'',
  Comm. Pure Appl. Math. 57, p. 657-703 (2004)
  
\bibitem{Ringstrom (2006)} H. Ringstr\"om,
  ``On curvature decay in expanding cosmological models'',
  Comm. Math. Phys. 264, p. 613-630 (2006)

\bibitem{Ringstrom (2013)}
  H. Ringstr\"om, \underline{On the topology and future stability of the
    universe}, Oxford University Press, Oxford (2013)

\bibitem{Ringstrom (2015)}
  H. Ringstr\"om, ``Instability of spatially homogeneous solutions in the
  class of $T^2$-symmetric solutions to Einstein's vacuum equations''
  Comm. Math. Phys. 334, p. 1299-1375 (2015)

    \end{thebibliography}
\end{document}